\documentclass[final,3p,times]{elsarticle}

\usepackage{graphicx}
\usepackage{amssymb}
\usepackage{amsthm}
\usepackage[fleqn]{amsmath}
\allowdisplaybreaks[1]
\biboptions{longnamesfirst,square,comma,compress}
\everymath{\displaystyle}
\usepackage{multicol}
\usepackage{array}
\usepackage{longtable} 

\begin{document}

\begin{frontmatter}

\title{Eighth-order Derivative-Free Family of Iterative Methods for Nonlinear Equations }
\author[rvt]{ Laila M Assas }
\ead{ lmassas@uqu.edu.sa}
\author[focal]{Fayyaz Ahmad \corref{cor1}\fnref{fn1} }
\ead{fayyaz.ahmad@upc.edu}
\author[els]{Malik Zaka Ullah }
\ead{ mzhussain@kau.edu.sa}

\cortext[cor1]{Corresponding author}
\fntext[fn1]{This research was supported by Spanish MICINN grants AYA2010-15685.}

\address[rvt]{Department of Mathematics, Umm al Qura University, Makkah, Kingdom of Saudi Arabia }
\address[focal]{Dept. de FÃ­sica i Enginyeria Nuclear,  Universitat PolitÃšcnica de Catalunya, Barcelona 08036, Spain }
\address[els]{ Department of Mathematics, King Abdulaziz University, Jeddah 21589, Kingdom of Saudi Arabia }

\begin{abstract}
In this note, we present an  eighth-order derivative-free family of iterative methods for nonlinear equations. The proposed family 
shows optimal  eight-order of convergence in the sense of the Kung and Traub conjecture \cite{5} and is based on the Steffensen 
derivative approximation used in the Newton-method. As a final step, having in mind computational purposes, a derivative-free polynomial 
base interpolation is used in order to get optimal order of convergence with only four functional evaluations. 
Numerical esperiments and few issues are discussed at the end of this note. 
\end{abstract}


\end{frontmatter}

\section{Introduction}
\label{Intro}
Let $f:D \subseteq \Re \longrightarrow  \Re $ be a sufficiently differentiable function of single variable in some neighborhood $D$ of
$ \alpha $, where $\alpha $ is a simple root ($f'(\alpha) \neq 0 $) of nonlinear algebraic equation $f(x)=0$. The well-known Newton method is denined by the iteration 
\begin{align} \label{eqns1}
 x_{n+1} &=  x_n - \frac{f(x_n)}{f'(x_n)}, 
\end{align}
which shows a second-order convergence. One can easily get Steffensen's approximation for first order derivative as 
\begin{align} \label{eqns2}
\begin{cases}
f(x_n-\kappa f(x_n)) &\approx  f(x_n)-\kappa f(x_n) f'(x_n), \\ 
\kappa f(x_n) f'(x_n)  &\approx  f(x_n)-  f(x_n-\kappa f(x_n)), \\
f'(x_n)  &\approx  \frac{1}{\kappa} \frac{f(x_n)-  f(x_n-\kappa f(x_n))}{ f(x_n)}.  
\end{cases}
\end{align}
If we substitute the derivative approximation (\ref{eqns2}) in (\ref{eqns1}), we obtain Steffensen's second order accurate derivative-free iterative method for non-linear equations \cite{1}. 
\begin{align}\label{eqns3}
\begin{cases}
w_n &= x_n-\kappa f(x_n),  \\
x_{n+1} &= x_n - \kappa  \frac{ f(x_n)^2}{f(x_n)-  f(w_n)}. 
\end{cases}
\end{align}

In 2012, an optimal eighth-order iterative method \cite{2} was proposed by Y. Khan et al. as 
\begin{align}\label{eqns4}
 \begin{cases}
  y_n &= x_n-\frac{f(x_n)}{f'(x_n)}, \\
   z_n &=  y_n -G\Bigg( \frac{f(y_n)}{f(x_n)} \Bigg) \frac{f(y_n)}{f'(x_n)}, \\
 x_{n+1} &= z_n - \frac{\mu}{\mu + \nu q_n^{2}} \frac{f(z_n)}{f'(z_n)},
 \end{cases}
\end{align}
where $\mu \neq  0$, $\nu \in  \Re $, $q_n = f(z_n)/f(x_n)$, $G(t)$ is a real-valued function with $G(0)=1$, $G'(0)=2$, $G''(0)< \infty$, and

\begin{align}  \label{eqns5}
\begin{cases}
 f'(z_n) &\approx K-C(y_n-z_n)-D(y_n-z_n)^2,\\ 
 H &= \frac{f(x_n)-f(y_n)}{x_n-y_n}, \\
 K &= \frac{f(x_n)-f(z_n)}{y_n-z_n}, \\
 D &= \frac{f^{'}(x_n)- H }{(x_n-y_n)(x_n-z_n)}-\frac{H-K}{(x_n-z_n)^2}, \\
 C &= \frac{H-K}{x_n-z_n}-D(   x_n+  y_n-2  z_n).
 \end{cases}
\end{align}
In the original draft of paper \cite{2} the expression for $C$ has typo-mistake, which is corrected here. Actually (\ref{eqns5}) 
polynomial interpolation  approximation for $f'(z_n)$ is given in  \cite{3}. Clearly (\ref{eqns4}) iterative 
scheme is not derivative free.
The main contribution in this paper is to use the idea of iterative scheme (\ref{eqns4}) by introducing Steffensen's derivative approximation 
for $f'(x_n)$ and then finally construct derivative-free approximation for $f'(z_n)$ without reducing order of convergence.
\section{Construction of derivative-free family}
First we construct an interpolation polynomial approximation for $f'(z_n)$. Suppose we have $f(x_n)$, $f(w_n)$ (defined in \ref{eqns3}), $f(y_n)$ 
and $f(z_n)$, One could construct a three-degree polynomial as follows
\begin{align} \label{eqns6}
 \begin{cases}
  p(\phi) &= f(y_n) + r_1 (\phi - y_n) + r_2  (\phi - y_n)^2 + r_3  (\phi - y_n)^3, \\
  p'(\phi) &=  r_1  + 2 r_2  (\phi - y_n) + 3 r_3  (\phi - y_n)^2. 
 \end{cases}
\end{align}
By using four functional values, we get the following system of equations: 
\begin{align} \label{eqns7}
 \begin{cases}
v_1 &=  r_1 a + r_2  a^2 + r_3  a^3, \\
v_2 &=  r_1 b + r_2  b^2 + r_3  b^3, \\
v_3 &=  r_1 c + r_2  c^2 + r_3  c^3,
 \end{cases}
\end{align}
where 
\begin{align}\label{eqns8}
 \begin{cases}
   v_1&=  f(x_n) - f(y_n ) ,\\
   v_2&=   f(z_n) - f(y_n ) , \\
   v_3&=  f(w_n) - f(y_n )  , \\
   a&=x_n-y_n , \\
   b&=z_n-y_n  \\
   c&=w_n-y_n .
 \end{cases}
\end{align}
 After solving (\ref{eqns7}) for $r_1$, $r_2$ and $r_3$ and substituting them in (\ref{eqns6}) implies the following approximation 
 for $f'(z_n)$:
\begin{align} \label{eqns9}
 f'(z_n) &\approx \psi_n =\frac{b (b-c) }{(a-b)(a-c)} \frac{v_1}{a}    +\frac{-3 b^2+2 b c+2 a b-a c }{(a-b)(b-c)} \frac{v_2}{b}+ 
 \frac{b (b-a)}{(a-c)(b-c)} \frac{v_3}{c}. 
 \end{align}

 We consider the following family of iterative methods for nonlinear equations:
 \begin{align} \label{eqns10}
  \begin{cases}
  w_n &= x_n-\kappa f(x_n),  \\
   y_n &=x_n - \kappa  \frac{ f(x_n)^2}{f(x_n)-  f(w_n)},\\
   z_n &= y_n - \kappa  \frac{ f(y_n) f(x_n)}{f(x_n)-  f(w_n)}   G(t_1,t_2),  \\
     x_{n+1} &= z_n - \frac{f(z_n)}{\psi_n}  H(s_1,s_2),  
  \end{cases}
 \end{align}
where  $t_1= \frac{f(y_n)}{f(x_n)}$, $t_2=\frac{f(y_n)}{f(w_n)}$, $s_1=\frac{f(z_n)}{f(x_n)}$, $s_2=\frac{f(z_n)}{f(w_n)}$ and 
$\kappa (\neq 0) \in \Re$.
\section{Convergence analysis}
We state the following theorem about the order of convergence of the family described in (\ref{eqns10}).

\newtheorem{mydef1}{Theorem}
\begin{mydef1}
 Let  $f: D\subseteq \Re \rightarrow \Re $ be a sufficiently differentiable function, and $\alpha \in D$ is a simple root of $f(x)=0$, 
 for an open interval $D$. If $x_0$ is chosen sufficiently close to $\alpha$, then the iterative scheme given in (\ref{eqns10}) converges to $\alpha$. If $G$ 
 and $H$ satisfy 
 \begin{align} \label{eqns11}
  G(0,0)&=1,\ \frac{\partial G}{\partial t_1} \Bigg|_{(0,0)} =1,\ \frac{\partial G}{\partial t_2} \Bigg|_{(0,0)} =1,\ 
  H(0,0)=1,\ \frac{\partial H}{\partial s_1}\Bigg|_{(0,0)} =0,\ \frac{\partial H}{\partial s_2} \Bigg|_{(0,0)} =0,   
 \end{align}
 and $\frac{\partial^2 G}{\partial t_1^2}$, $ \frac{\partial^2 G}{\partial t_2^2}$ ,$ \frac{\partial^2 G}{ \partial t_1  \partial t_2 }$, 
 $\frac{\partial^2 H}{\partial s_1^2}$, $ \frac{\partial^2 H}{\partial s_2^2}$ ,$ \frac{\partial^2 H}{ \partial s_1  \partial s_2 }$
are bounded at $(0,0)$
then the iterative scheme (\ref{eqns10}) shows an order of  convergence at least equal to eight.
\end{mydef1}
\begin{proof}
Let the error at step $n$ be denoted by $e_n= x_n-\alpha$ and let us define $c_1 = f'(\alpha)$ and $c_k = \frac{1}{k!} \frac{f^{(k)}(\alpha)}{f'(\alpha)}$, 
$k = 2,3,\cdots$. If we expand $f$ around the root $\alpha$ and express it in terms of powers of error $e_n$, we obtain
\begin{align}  
 f(x_n) &= c_1 ( e_n+ c_2  e_n^2+ c_3  e_n^3+ c_4  e_n^4+ c_5  e_n^5+ c_6  e_n^6+ c_7  e_n^7+ c_8  e_n^8+ O( e_n^9) ),  \label{eqns12}  \\ 
 \begin{split}  \label{eqns13} 
 f(w_n) &= -c_1 (-1+\kappa c_1) e_n + c_1 c_2 (-3  \kappa c_1+1+ \kappa^2 c_1^2) e_n^2 -c_1 (4  \kappa c_1 c_3+2 c_2^2  \kappa c_1-2  \kappa^2 c_1^2 c_2^2-c_3-3 c_3  \kappa^2 c_1^2 \\
 &\quad +c_3  \kappa^3 c_1^3) e_n^3 +  c_1 (-5 \kappa  c_1  c_4-5  c_2 \kappa  c_1  c_3+8 \kappa^2  c_1^2  c_2  c_3+\kappa^2  c_1^2  c_2^3-3  c_2 \kappa^3  c_1^3  c_3+ c_4+6  c_4 \kappa^2  c_1^2-4  c_4 \kappa^3  c_1^3 \\
 &\quad + c_4 \kappa^4  c_1^4) e_n^4 +\cdots + O(e_n^9),  
 \end{split} \\
 \begin{split} \label{eqns14} 
  y_n -\alpha &= - c_2  (-1+  \kappa   c_1) e_n^2 + (2  c_3-3   \kappa  c_1  c_3+ c_3   \kappa^2  c_1^2+2  c_2^2   \kappa  c_1-2  c_2^2-  \kappa^2  c_1^2  c_2^2 ) e_n^3  + (3  c_4+10  c_2  \kappa  c_1  c_3-6  \kappa  c_1  c_4 \\
   &\quad +4  c_4  \kappa^2  c_1^2- c_4  \kappa^3  c_1^3-7  \kappa^2  c_1^2  c_2  c_3-7  c_2  c_3-5  \kappa  c_1  c_2^3+2  c_2  \kappa^3  c_1^3  c_3+3  \kappa^2  c_1^2  c_2^3+4  c_2^3- \kappa^3  c_1^3  c_2^3) e_n^4 \\
   &\quad +\cdots+O(e_n^9),
 \end{split} \\
 \begin{split} \label{eqns15}
  f(y_n) &= - c_1  c_2 (-1+ \kappa  c_1) e_n^2 -  c_1  (-2    c_3+3   \kappa    c_1    c_3-  c_3   \kappa^2    c_1^2-2
  c_2^2   \kappa    c_1+2    c_2^2+ \kappa^2    c_1^2    c_2^2) e_n^3 -c_1 (-3 c_4-10 c_2 \kappa c_1 c_3\\
  &\quad +6 \kappa c_1 c_4-4 c_4 \kappa^2 c_1^2+c_4 \kappa^3 c_1^3+7 \kappa^2 c_1^2 c_2 c_3+7 c_2 c_3+7 \kappa c_1 c_2^3-2 c_2 \kappa^3 c_1^3 c_3-4 \kappa^2 c_1^2 c_2^3-5 c_2^3 \\
  &\quad +\kappa^3 c_1^3 c_2^3) e_n^4 + \cdots + O(e_n^9),
 \end{split} \\
 \begin{split} \label{eqns16}
  \frac{f(y_n)}{f(x_n)} &=  -c_2 (-1+ \kappa c_1) e_n + ( 
  2  c_3-3  \kappa  c_1  c_3+ c_3  \kappa^2  c_1^2+3  c_2^2  \kappa  c_1-3  c_2^2- \kappa^2  c_1^2  c_2^2 ) e_n^2 +\cdots+O(e_n^9),
 \end{split} \\
 \begin{split} \label{eqns17}
  \frac{f(y_n)}{f(w_n)} &= c_2 e_n + (-  \kappa c_1 c_3+2 c_2^2   \kappa c_1+2 c_3-3 c_2^2) e_n^2 + \cdots+O(e_n^9).
 \end{split}
\end{align}
The Taylor series expansion of $G(t_1,t_2)$ is given by 
\begin{align} \label{eqns18}
 \begin{split}
 G\Bigg(\frac{f(y_n)}{f(x_n)},\frac{f(y_n)}{f(w_n)}\Bigg) &= 1 +  \frac{f(y_n)}{f(x_n)} +
 \frac{f(y_n)}{f(w_n)} + A_1 \Bigg( \frac{f(y_n)}{f(x_n)}  \Bigg)^2 + A_2 \Bigg( \frac{f(y_n)}{f(w_n)} 
 \Bigg)^2 + A_3 \Bigg( \frac{f(y_n)}{f(x_n)}  \Bigg)  \Bigg( \frac{f(y_n)}{f(w_n)}  \Bigg) + O\Big( t_1^3,t_2^3 \Big).   
 \end{split}
\end{align}
By using (\ref{eqns12}), (\ref{eqns13}), (\ref{eqns15}), (\ref{eqns16}), (\ref{eqns17}), we find

\begin{align}
%  \begin{split} 
  z_n-\alpha &= (- A_1 c_2^2+6   \kappa^2 c_1^2 c_2^2-c_3+2   \kappa c_1 c_3-10 c_2^2   \kappa c_1-c_3   \kappa^2 c_1^2+5 c_2^2+ A_1 c_2^2   \kappa c_1+2 c_2^2  A_3   \kappa c_1-c_2^2  A_3   \kappa^2 c_1^2\nonumber \\
   &\quad +3  A_2 c_2^2   \kappa c_1-3  A_2 c_2^2   \kappa^2 c_1^2+ A_2 c_2^2   \kappa^3 c_1^3- A_2 c_2^2-  \kappa^3 c_1^3 c_2^2-c_2^2  A_3) c_2 e_n^4 + (-4  \kappa^2  c_1^2  c_3^2-4  \kappa^2  c_1^2  c_4  c_2\nonumber \\
   &\quad -31  A_2  c_2^4  \kappa  c_1+36  A_2  c_2^4  \kappa^2  c_1^2-19  A_2  c_2^4  \kappa^3  c_1^3+4  A_2  c_2^4  \kappa^4  c_1^4-23  c_2^4  A_3  \kappa  c_1+18  c_2^4  A_3  \kappa^2  c_1^2-5  c_2^4  A_3  \kappa^3  c_1^3\nonumber \\  
   &\quad -15  A_1  c_2^4  \kappa  c_1+6  A_1  c_2^4  \kappa^2  c_1^2+5  c_4  \kappa  c_1  c_2+3  \kappa^4  c_1^4  c_2^2  c_3+ c_4  \kappa^3  c_1^3  c_2+68  \kappa^2  c_1^2  c_2^2  c_3-25  c_3  \kappa^3  c_1^3  c_2^2-2  c_4  c_2\nonumber \\
   &\quad -78  \kappa  c_1  c_2^2  c_3-3  \kappa^4  c_1^4  c_2^4+5  \kappa  c_1  c_3^2+ c_3^2  \kappa^3  c_1^3-2  c_3^2-36  c_2^4+32  c_2^2  c_3-66  c_2^4  \kappa^2  c_1^2+80  c_2^4  \kappa  c_1+24  \kappa^3  c_1^3  c_2^4\nonumber \\
   &\quad +9  A_1  c_2^2  \kappa  c_1  c_3-3  A_1  c_2^2  c_3  \kappa^2  c_1^2+21  A_2  c_2^2  \kappa  c_1  c_3-27  A_2  c_2^2  c_3  \kappa^2  c_1^2 +15  A_2  c_2^2  c_3  \kappa^3  c_1^3-3  A_2  c_2^2  \kappa^4  c_1^4  c_3 \nonumber \\
   &\quad +15  c_2^2  A_3  \kappa  c_1  c_3-12  c_2^2  A_3  c_3  \kappa^2  c_1^2+3  c_2^2  A_3  c_3  \kappa^3  c_1^3-6  A_1  c_2^2  c_3-6  A_2  c_2^2  c_3-6  c_2^2  A_3  c_3 +10  A_1  c_2^4+10  c_2^4  A_3 \nonumber \\
   &\quad +10  A_2  c_2^4) e_n^5 + \cdots + O(e_n^9), \label{eqns19}
%  \end{split} \\
\end{align}

\begin{align}
 \begin{split} \label{eqns20}
  f(z_n) &=   c_1 (-A_1  c_2^2+6  \kappa^2  c_1^2  c_2^2- c_3+2  \kappa  c_1  c_3-10  c_2^2  \kappa  c_1- c_3  \kappa^2  c_1^2+5  c_2^2+A_1  c_2^2  \kappa  c_1+2  c_2^2 A_3  \kappa  c_1- c_2^2 A_3  \kappa^2  c_1^2\\
  &\quad  +3 A_2  c_2^2  \kappa  c_1-3 A_2  c_2^2  \kappa^2  c_1^2+A_2  c_2^2  \kappa^3  c_1^3-A_2  c_2^2- \kappa^3  c_1^3  c_2^2- c_2^2 A_3)  c_2 e_n^4 + c_1 (-4 \kappa^2 c_1^2 c_3^2-4 \kappa^2 c_1^2 c_4 c_2 \\
  &\quad  -31 A_2 c_2^4 \kappa c_1+36 A_2 c_2^4 \kappa^2 c_1^2-19 A_2 c_2^4 \kappa^3 c_1^3+4 A_2 c_2^4 \kappa^4 c_1^4-23 c_2^4 A_3 \kappa c_1+18 c_2^4 A_3 \kappa^2 c_1^2-5 c_2^4 A_3 \kappa^3 c_1^3\\
  &\quad  -15 A_1 c_2^4 \kappa c_1+6 A_1 c_2^4 \kappa^2 c_1^2+5 c_4 \kappa c_1 c_2+3 \kappa^4 c_1^4 c_2^2 c_3+c_4 \kappa^3 c_1^3 c_2+68 \kappa^2 c_1^2 c_2^2 c_3-25 c_3 \kappa^3 c_1^3 c_2^2-2 c_4 c_2\\
  &\quad  -78 \kappa c_1 c_2^2 c_3-3 \kappa^4 c_1^4 c_2^4+5 \kappa c_1 c_3^2+c_3^2 \kappa^3 c_1^3-2 c_3^2-36 c_2^4+32 c_2^2 c_3-66 c_2^4 \kappa^2 c_1^2+80 c_2^4 \kappa c_1+24 \kappa^3 c_1^3 c_2^4\\
  &\quad  +9 A_1 c_2^2 \kappa c_1 c_3-3 A_1 c_2^2 c_3 \kappa^2 c_1^2+21 A_2 c_2^2 \kappa c_1 c_3-27 A_2 c_2^2 c_3 \kappa^2 c_1^2+15 A_2 c_2^2 c_3 \kappa^3 c_1^3-3 A_2 c_2^2 \kappa^4 c_1^4 c_3 \\
  &\quad +15 c_2^2 A_3 \kappa c_1 c_3-12 c_2^2 A_3 c_3 \kappa^2 c_1^2+3 c_2^2 A_3 c_3 \kappa^3 c_1^3-6 A_1 c_2^2 c_3-6 A_2 c_2^2 c_3-6 c_2^2 A_3 c_3\\
  &\quad  +10 A_1 c_2^4+10 c_2^4 A_3+10 A_2 c_2^4) e_n^5 + \cdots +O(e_n^9), 
 \end{split}\\
 \begin{split} \label{eqns21}
  \psi_n &=  (2 A_2 c_2^3 \kappa^3 c_1^3-2 \kappa^3 c_1^3 c_2^3+12 \kappa^2 c_1^2 c_2^3-2 \kappa^2 c_1^2 c_2 c_3-6 A_2 c_2^3 \kappa^2 c_1^2+c_4 \kappa^2 c_1^2-2 c_2^3 A_3 \kappa^2 c_1^2-20 \kappa c_1 c_2^3\\
  &\quad +6 A_2 c_2^3 \kappa c_1-2 \kappa c_1 c_4+2 A_1 c_2^3 \kappa c_1+4 c_2 \kappa c_1 c_3+4 c_2^3 A_3 \kappa c_1+10 c_2^3-2 A_1 c_2^3 -2 c_2^3 A_3-2 c_2 c_3 \\
  &\quad -2 A_2 c_2^3+c_4) c_1 c_2  e_n^4 + \cdots + O(e_n^9).
 \end{split}
\end{align} 
Finally, $H(s_1,s_2)$ has Taylor's expansion 
\begin{align} 
 \begin{split}\label{eqns22}
  H\Bigg(\frac{f(z_n)}{f(x_n)},\frac{f(z_n)}{f(w_n)}\Bigg) &= 1+ B_1 
  \Bigg( \frac{f(z_n)}{f(x_n)} \Bigg)^2 + B_2 \Bigg( \frac{f(z_n)}{f(w_n)} \Bigg)^2 + 
  B_3  \Bigg( \frac{f(z_n)}{f(x_n)} \Bigg) \Bigg( \frac{f(z_n)}{f(w_n)} \Bigg)+O(s_1^3,s_2^3).
 \end{split}
\end{align}
From (\ref{eqns20}) and (\ref{eqns21}), we deduced the following error equation which leads to the desired result
\begin{align} 
%  \begin{split} 
  e_{n+1} &=  c_2^2 (-6  \kappa^2  c_1^2  c_3  c_4-10 A_1  c_2^5-10 A_2  c_2^5-10  c_2^5 A_3+6  c_2  \kappa^2  c_1^2  c_3^2+31  \kappa^2  c_1^2  c_4  c_2^2-4  c_2  \kappa^3  c_1^3  c_3^2\nonumber \\
  &\quad +4  c_3  \kappa  c_1  c_4+4  c_3  c_4  \kappa^3  c_1^3+46  c_3  \kappa^3  c_1^3  c_2^3-23  c_4  \kappa^3  c_1^3  c_2^2+8  c_4  \kappa^4  c_1^4  c_2^2- c_4  \kappa^4  c_1^4  c_3-4  c_2  \kappa  c_1  c_3^2-20  \kappa  c_1  c_4  c_2^2\nonumber \\
  &\quad - c_3  c_4-62  c_3  \kappa^2  c_1^2  c_2^3+ \kappa^4  c_1^4  c_2  c_3^2-16  \kappa^4  c_1^4  c_2^3  c_3-12  \kappa^5  c_1^5  c_2^5+2  \kappa^5  c_1^5  c_3  c_2^3- \kappa^5  c_1^5  c_4  c_2^2+6 A_1  c_2^3  \kappa^2  c_1^2  c_3\nonumber \\
  &\quad +3 A_1  c_2^2  \kappa  c_1  c_4-6 A_1  c_2^3  \kappa  c_1  c_3-2 A_1  c_2^3  \kappa^3  c_1^3  c_3-3 A_1  c_2^2  c_4  \kappa^2  c_1^2+A_1  c_2^2  c_4  \kappa^3  c_1^3-8  c_2^3 A_3  \kappa^3  c_1^3  c_3+4  c_2^2 A_3  \kappa  c_1  c_4\nonumber \\
  &\quad -6  c_2^2 A_3  \kappa^2  c_1^2  c_4+12  c_2^3 A_3  \kappa^2  c_1^2  c_3+40  c_2^3  c_3  \kappa  c_1-A_1  c_2^2  c_4+2 A_1  c_2^3  c_3-A_2  c_2^2  c_4+2 A_2  c_2^3  c_3- c_2^2 A_3  c_4\nonumber \\
  &\quad +2  c_2^3 A_3  c_3+25  c_2^5+2  c_2^3 A_3  \kappa^4  c_1^4  c_3+4  c_2^2 A_3  \kappa^3  c_1^3  c_4- c_2^2 A_3  \kappa^4  c_1^4  c_4-8  c_2^3 A_3  \kappa  c_1  c_3-20 A_2  c_2^3  \kappa^3  c_1^3  c_3\nonumber \\
  &\quad +5 A_2  c_2^2  \kappa  c_1  c_4-10 A_2  c_2^2  \kappa^2  c_1^2  c_4+20 A_2  c_2^3  \kappa^2  c_1^2  c_3+10 A_2  c_2^3  \kappa^4  c_1^4  c_3+10 A_2  c_2^2  \kappa^3  c_1^3  c_4-5 A_2  c_2^2  \kappa^4  c_1^4  c_4  \nonumber \\ 
  &\quad-10 A_2  c_2^3  \kappa  c_1  c_3-2 A_2  c_2^3  \kappa^5  c_1^5  c_3+A_2  c_2^2  \kappa^5  c_1^5  c_4+ c_2  c_3^2-10  c_2^3  c_3+5  c_2^2  c_4+160  \kappa^2  c_1^2  c_2^5-130  \kappa^3  c_1^3  c_2^5\nonumber \\
  &\quad -100  c_2^5  \kappa  c_1+56  \kappa^4  c_1^4  c_2^5-32 A_1  c_2^5  \kappa^2  c_1^2+30 A_1  c_2^5  \kappa  c_1+14 A_1  c_2^5  \kappa^3  c_1^3+46  c_2^5 A_3  \kappa^3  c_1^3-62  c_2^5 A_3  \kappa^2  c_1^2\nonumber \\
  &\quad -16  c_2^5 A_3  \kappa^4  c_1^4+40  c_2^5 A_3  \kappa  c_1+108 A_2  c_2^5  \kappa^3  c_1^3-102 A_2  c_2^5  \kappa^2  c_1^2-62 A_2  c_2^5  \kappa^4  c_1^4+50 A_2  c_2^5  \kappa  c_1+18 A_2  c_2^5  \kappa^5  c_1^5\nonumber \\
  &\quad +A_1^2  c_2^5+A_2^2  c_2^5+ c_2^5 A_3^2-2 A_1^2  c_2^5  \kappa  c_1-2 A_1  c_2^5  \kappa^4  c_1^4+A_1^2  c_2^5  \kappa^2  c_1^2-4  c_2^5 A_3^2  \kappa^3  c_1^3+6  c_2^5 A_3^2  \kappa^2  c_1^2-4  c_2^5 A_3^2  \kappa  c_1\nonumber \\
  &\quad +2  c_2^5 A_3  \kappa^5  c_1^5+ c_2^5 A_3^2  \kappa^4  c_1^4+15 A_2^2  c_2^5  \kappa^4  c_1^4-20 A_2^2  c_2^5  \kappa^3  c_1^3+15 A_2^2  c_2^5  \kappa^2  c_1^2-6 A_2^2  c_2^5  \kappa  c_1-6 A_2^2  c_2^5  \kappa^5  c_1^5 \nonumber \\
  &\quad +A_2^2  c_2^5  \kappa^6  c_1^6-2 A_2  c_2^5  \kappa^6  c_1^6-8 A_1  c_2^5 A_2  \kappa^3  c_1^3+12 A_1  c_2^5 A_2  \kappa^2  c_1^2+6 A_1  c_2^5 A_3  \kappa^2  c_1^2-8 A_1  c_2^5 A_2  \kappa  c_1-6 A_1  c_2^5 A_3  \kappa  c_1 \nonumber \\
  &\quad +2 A_1  c_2^5  \kappa^4  c_1^4 A_2-2 A_1  c_2^5  \kappa^3  c_1^3 A_3+10  c_2^5 A_3  \kappa^4  c_1^4 A_2-20  c_2^5 A_3  \kappa^3  c_1^3 A_2+20  c_2^5 A_3  \kappa^2  c_1^2 A_2-10  c_2^5 A_3  \kappa  c_1 A_2 \nonumber \\
  &\quad -2  c_2^5 A_3  \kappa^5  c_1^5 A_2+2 A_1  c_2^5 A_3+2 A_1  
  c_2^5 A_2+2 A_2  c_2^5 A_3+ \kappa^6  c_1^6  c_2^5) e_n^8 + O(e_n^9). \label{eqns23}
%  \end{split} 
\end{align}
\end{proof}

It is clear that the considered family of numerical schemes requires four functional evaluations and attains optimal convergence order eight according to  Kung and Traub conjecture
which can be stated as follows \cite{5}:  if $n$ is the total number of functional evaluations per iteration, then the optimal convergence order of the associated numerical procedure is $2^{n-1}$. 

\section{Numerical Results}
\newtheorem{mydef3}{Definition}
\begin{mydef3}
 The computational order of convergence \cite{4}, can be approximated by
 \begin{align} \label{eqns24}
 COC &\approx  \frac{ln|(x_{n+1}-\alpha)(x_n-\alpha)^{-1}|}{ln|(x_n-\alpha)(x_{n-1}-\alpha)^{-1}|},
 \end{align}
where $x_{n-1}$, $x_n$ and $x_{n+1}$ are successive iterations closer to the root $\alpha$ of $f(x)=0$.
\end{mydef3}
For the purpose of comparison between newly developed family and other derivative-free methods, a list of derivative-free methods for nonlinear 
equations is presented here. 

\subsection{The Kung-Traub Eighth-order Derivative-free Method (K-T)}
The Kung-Traub eighth-order derivative-free  method is discussed in \cite{5,6}, and also considered in \cite{7} is given as 

 \begin{align}\label{eqns25}
 \begin{split}
 \begin{cases}
  w_n &= x_n + \beta f(x_n), \\
  y_n &= x_n -\Bigg( \frac{\beta f(x_n)^2 }{f(w_n)-f(x_n)}   \Bigg), \\
  z_n &= y_n - \Bigg(   \frac{f(x_n)f(w_n)}{f(y_n)-f(x_n)}  \Bigg) \Bigg[ \frac{1}{f[w_n,x_n]} -\frac{1}{f[w_n,x_n]}  \Bigg],\\
  x_{n+1} &= z_n - \Bigg( \frac{f(w_n)f(x_n)f(y_n)}{f(z_n)-f(x_n)}\Bigg)  \\
  &\quad \Bigg\{ \Bigg( \frac{1}{f(z_n)-f(w_n)}  \Bigg) \Bigg[ \frac{1}{f[y_n,z_n]} -\frac{1}{f[w_n,y_n]} 
  \Bigg]   - \Bigg( \frac{1}{f(y_n)-f(x_n)} \Bigg)   \Bigg[ \frac{1}{f[w_n,y_n]} -\frac{1}{f[w_n,x_n]}  \Bigg]    \Bigg\}.
 \end{cases}
 \end{split}
 \end{align}

\subsection{R. Thukral  $M1$, $M2$, $M3$ Methods}
In 2011, R. Thukral \cite{7} presented three variants of his proposed eighth-order three-point derivative-free method. Three members
of the family called by author namely, $M1$, $M2$, and $M3$, are listed as 
\begin{align} 
 \phi_1 &= \Bigg(  1-\frac{f(y_n)}{f(w_n)}  \Bigg)^{-1}, \label{eqns26} \\
 \phi_2 &= \Bigg( 1+\frac{f(y_n)}{f(w_n)} + \Big( \frac{f(y_n)}{f(w_n)}  \Big)^2  \Bigg), \label{eqns27} \\
 \phi_3 &=   \frac{f[x_n,w_n]}{f[w_n,y_n]}, \label{eqns28} 
\end{align}
and 
\begin{align} \label{eqns29}
\begin{cases}
 w_n &= x_n +\beta f(x_n), \\
 y_n &= x_n-\Bigg( \frac{\beta f(x_n)^2}{f(w_n)-f(x_n)}\Bigg), \\
 z_n &= y_n - \phi_k \Bigg( \frac{f(y_n)}{f[x_n,y_n]} \Bigg), \\
 x_{n+1}&= z_n - \Bigg(  1-\frac{f(z_n)}{f(w_n)} \Bigg)^{-1} \Bigg(  
 1- \frac{f(y_n)^3}{f(w_n)^2 f(x_n)} \Bigg) \Bigg( \frac{f[x_n,y_n] f(z_n)}{f[y_n,z_n] f[x_n,z_n]}\Bigg),
 \end{cases}
\end{align}

\begin{table}[!htbp]
\begin{center}
 \begin{tabular}{p{7cm}  p{5cm}}
  \hline 
  \hline \\
  Functions & Roots \\
  \hline 
  \hline \\
$f_1 (x) = exp(x)\ sin(x) + ln(1 + x^2 )$    		&  $\alpha=0 $ \\ 
$f_2 (x) = x^{15} + x^4 + 4x^2 -15  		$	&  $\alpha = 1.148538 . . .$  \\
$f_3 (x) = (x - 2)(x^{10} + x + 1)\ exp(-x-1)$		&  $\alpha=2  $\\
$f_4 (x) = exp(-x^2 + x + 2) - cos(x + 1) + x^3 + 1$	&  $\alpha = -1$  \\
$f_5 (x) = (x+1)\ exp(sin(x))-x^2 \ exp(cos(x))-1$	&  $\alpha=0  $\\
$f_6 (x) = sin(x)^2 - x^2 + 1		$		&  $\alpha = 1.40449165 . . . $ \\
$f_7 (x) = 10\ exp(-x^2 ) -1		$		&  $\alpha = 1.517427 . . .  $\\
$f_8 (x) = (x^2 -1)^{-1} - 1		$		&  $\alpha = 1.414214 . . .  $\\
$f_{9} (x) = ln(x^2 + x + 2) - x + 1	$		&  $\alpha = 4.15259074 . . .$  \\
$f_{10} (x) = cos(x)^2 - x/5		$	&  $\alpha = 1.08598268 . . .$  \\
$f_{11} (x) = sin(x) - \frac{x}{2}		$       &  $\alpha=0 $ \\
$f_{12} (x) = x^{10} - 2x^3 - x + 1		$	&  $\alpha = 0.591448093 . . .$  \\
$f_{13} (x) = exp(sin(x)) - x + 1         $            &  $\alpha = 2.63066415 . . . $\\ 
\end{tabular}
 \begin{tabular}{p{1.5cm}  p{2.5cm}  p{1.5cm} p{1.5cm} p{1.5cm}  p{1.5cm} p{1.5cm}  p{1.5cm} }
  \hline 
  \hline \\
   ($f_n(x)$,$x_0$ ) & L  & K-T & M1 & M2 & M3 & P1 & P2  \\   
  \hline 
  \hline \\
   $f_1$,\ 0.25  & (L1) 6.38e-247    & 3.14e-136 & 1.69e-141 & 7.43e-142 & 1.69e-141 & 3.20e-113 & 8.98e-120 \\ 
   $f_2$,\ 1.1   & (L1) 1.2376e-652  & 3.72e-61 & 3.44e-62 & 3.44e-62 & 3.44e-62 & 2.68e-7 & 2.67e-7 \\ 
   $f_3$,\ 2.1   & (L1) 1.057e-422   & 1.91e-60 & 1.49e-60 & 1.49e-60 & 1.49e-60 & 7.71e-8 & 7.56e-8 \\ 
   $f_4 $,\ -0.5  & (L1) 2.952e-383   & 5.11e-362 & 1.92e-362 & 1.93e-362 & 1.92e-362 & 9.99e-367 & 8.78e-366 \\ 
   $f_5 $,\ 0.25  & (L1) 2.336e-407   & 4.13e-328 & 6.52e-326 & 9.47e-326 & 6.52e-326 & 1.98e-322 & 2.56e-332 \\ 
   $f_6 $,\ 1.2   & (L8) 1.719e-421   & 1.00e-327 & 4.58e-341 & 7.57e-344 & 4.58e-341 & 1.79e-381 & 1.72e-405 \\ 
   $f_7 $,\ 2     & (L2) 7.264e-238   & 5.19e-88 & 1.24e-120 & 6.40e-124 & 1.24e-120 & 1.51e-187 & 6.79e-228 \\ 
   $f_8 $,\ 1.7     &  (L3) 1.429e-234    & 1.23e-113 & 1.74e-171 & 5.45e-188 & 1.74e-171 & 5.96e-211 & 4.84e-167 \\ 
   $f_{9} $,\ 4.4  &  (L4) 2.504e-997    & 1.15e-928 & 4.52e-942 & 1.27e-965 & 4.52e-941 & 6.15e-904 & 4.11e-937 \\ 
   $f_{10} $,\ 1.5  &  (L5) 2.81e-305    & 7.19e-303 & 5.07e-284 & 1.84e-245 & 5.07e-285 & 4.91e244 & 1.78e-275 \\ 
   $f_{11} $,\ 0.25  &  (L6) 2.35e-1143    & 3.65e-782 & 1.00e-819 & 4.98e-823 & 1.00e-819 & 5.13e-794 & 5.13e-812 \\ 
   $f_{12}$,\ 0.25  &  (L6) 7.86e-318    & 2.03e-256 & 5.65e256 & 1.82e-254 & 5.65e-256 & 1.07e-264 & 6.31e-268 \\ 
   $f_{13}$,\ 2.0  &  (L7) 2.54e-436    & 2.63e-396 & 1.94e-378 & 5.1e-378 & 1.94e-378 & 8.70e-380 & 6.80e-379 \\ 
\hline  
\hline  \\
    &      &         &       (COC)         &      &   &   &     \\ 
\hline
\hline  
  $f_1$     & (L1) 7.9999   & 7.9986 & 7.9995 & 7.9998  &7.9995	& 7.9958	& 7.9978\\ 
   $f_2$    & (L1) 8.0000   & 7.8671 & 7.9371 & 7.9371  &7.9371 	& 3.2715	& 3.2731\\ 
   $f_3$    & (L1) 7.9999   & 7.8660 & 7.9047 & 7.9047  &7.9047 	&4.2595 	& 4.2675\\ 
   $f_4 $   & (L1) 8.0000   & 7.9905 & 7.9905 & 7.9905  &7.9905 	&7.9907 	& 7.9907\\ 
   $f_5 $   & (L1) 8.0000   & 7.9884 & 7.9882 & 7.9882  &7.9882 	&8.0000 	& 8.0000\\ 
   $f_6 $   & (L8) 8.0000   & 8.000 &  8.0000&  8.0000  &8.0000 	&8.0000 	&8.0000 \\ 
   $f_7 $   & (L2) 7.9999   & 8.0097 & 8.0025 & 8.0018  &8.0025 	&8.0005 	&8.0001 \\ 
   $f_8 $   & (L3) 8.0000   & 8.0027 & 8.0004 & 8.0002  &8.0004 	&8.0000 	&8.0002 \\ 
   $f_{9}$ & (L4) 8.0000   & 8.0000 & 8.0000 & 8.0000  &8.0000 	&8.0000 	&8.0000 \\ 
   $f_{10}$ & (L5) 8.0000   & 8.0002 & 7.9845 & 7.9797  &7.9845 	&8.0003 	&7.9833 \\ 
   $f_{11}$ & (L6) 11.000   & 10.996 & 10.996 & 10.996  &10.996 	&10.996 	&10.996 \\ 
   $f_{12}$ & (L6) 8.0000   & 8.0000 & 7.9809 & 7.9807  &7.9809 	&7.9822 	&8.0000 \\ 
   $f_{13}$ & (L7) 7.9999   & 8.0000 & 8.0000 & 8.0000  &8.0000 	&8.0000 	&8.0000 \\ 
   \hline
\end{tabular}
\end{center}
\caption{Numerical comparison between three-point derivative-free methods }\label{tb2:numerical comparison between methods}
\end{table}

where $k=1,2,3$, $\beta \in \Re^{+}$, $\phi_k$ are listed in (26)-(28). (\ref{eqns29}) is called $M1$, $M2$ and $M3$ for $\phi_1$, $\phi_2$ and
$\phi_3$ respectively.
\subsection{Petkovic et al. Type Methods}
In \cite{7}, author developed Petkovic type 1 (P1) and type 2 (P2) derivative-free methods for the comparison of numerical efficiency, (P1)
and (P2) respectively, are written as 
\begin{align}\label{eqns30}
 \begin{cases}
 \begin{split}
 w_n &= x_n +\beta f(x_n), \\
 y_n &= x_n-\Bigg( \frac{\beta f(x_n)^2}{f(w_n)-f(x_n)}\Bigg), \\
 z_n &= y_n-\Bigg( 1+ \frac{f(y_n)}{f(w_n)}+\frac{f(y_n)}{f(x_n)}  \Bigg) \Bigg[ \frac{(w_n-x_n)f(y_n)}{f(w_n)-f(x_n)}  \Bigg], \\
 x_{n+1} &= z_n - \Bigg( 1-\frac{f(z_n)}{f(w_n)}  \Bigg)^{-1} \\
 &\quad \Bigg( 1- \frac{2 f(y_n)^3 }{f(w_n)^2 f(x_n)} - \frac{f(y_n)^3}{f(w_n)f(x_n)^2} - \Bigg( \frac{f(y_n)}{f(w_n)} \Bigg)^3  \Bigg) 
 \Bigg( \frac{f[x_n,y_n]f(z_n)}{f[y_n,z_n]f[x_n,z_n]} \Bigg),
 \end{split}
 \end{cases}
\end{align}
and
\begin{align}\label{eqns31}
 \begin{cases}
 \begin{split}
 w_n &= x_n +\beta f(x_n), \\
 y_n &= x_n-\Bigg( \frac{\beta f(x_n)^2}{f(w_n)-f(x_n)}\Bigg), \\
 z_n &= y_n-\Bigg( \frac{1+f(y_n)f(x_n)^{-1}}{1-f(y_n)f(w_n)^{-1}}   \Bigg) \Bigg( \frac{f(y_n) (w_n-x_n)}{f(w_n)-f(x_n)} \Bigg), \\
 x_{n+1} &= z_n - \Bigg( 1-\frac{f(z_n)}{f(w_n)}  \Bigg)^{-1}  \Bigg( 1- \frac{2 f(y_n)^3 }{f(w_n)^2 f(x_n)} - \frac{f(y_n)^3}{f(w_n)f(x_n)^2}  \Bigg) 
 \Bigg( \frac{f[x_n,y_n]f(z_n)}{f[y_n,z_n]f[x_n,z_n]} \Bigg).
 \end{split}
 \end{cases}
\end{align}
\subsection{Proposed family (L)}
We define the following weight functions:
\begin{align} 
G_1(t_1,t_2)    &= \frac{1}{Ã1-(t_1+t_2)+\omega (t_1+t_2)^2}, \ \omega \in \Re,  \label{eqns32} \\
G_2 (t_1,t_2) &= 1 + t_1 + t_2 +  t_1^2 + 1.9 t_2^2 + 4.4 t_1 t_2,  \label{eqns33} \\
H_1(s_1,s_2) &= 1, \label{eqns34} \\
H_2(s_1,s_2) &= \frac{1}{1 + s_1 s_2 + s_1^2 + s_2^2},   \label{eqns35} \\
H_3(s_1,s_2) &= 1 + s_2^4 + s_2^6, \label{eqns36} \\
H_4(s_1,s_2) &= 1 + s_1^2 + s_2^2 + 2 s_1 s_2, \label{eqns37} \\
H_5(s_1,s_2) &= \frac{1}{1-20 s_1 s_2}, \label{eqns39} 
\end{align}
where $t_i$ and $s_i$ are defined in (\ref{eqns10}). Further we give names to methods for the purpose of simplicity as follows
\begin{align}\label{eqns38}
\begin{cases}
L1 = (G_1,\ H_1,\ \omega=+0.01,\ \kappa=0.01),\  L2 = (G_1,\ H_1,\ \omega=-0.022,\ \kappa=0.01), \\ 
L3 = (G_1,\ H_1,\ \omega=-0.001,\ \kappa=0.01),\  L4 = (G_2,\ H_1,\ \omega=+0.01,\ \kappa=0.01), \\ 
L5 = (G_1,\ H_3,\ \omega=-0.01,\ \kappa=0.01),\  L6 = (G_1,\ H_2,\ \omega=+0.01,\ \kappa=0.01), \\
L7 = (G_1,\ H_4,\ \omega=+0.01,\ \kappa=0.01),\  L8 = (G_1,\ H_5,\ \omega=+0.01,\ \kappa=0.01).
\end{cases}
\end{align}

A set of thirteen nonlinear equations is used for numerical computations from \cite{7}, 
in Table \ref{tb2:numerical comparison between methods}. All the families in the numerical implementation
are derivative-free and use four function evaluations to get the order of convergence eight . 
For all methods, 12 (TNFE) total number of function evaluations are used, and absolute error ($|x_n-\alpha|$) is displayed. Computational 
order of convergence is calculated according to (\ref{eqns24}) for the method. All numerical values for 
methods K-T, M1, M2, M3, P1, P2 are taken from \cite{7}.

\section{Conclusion}
In this note, we have presented a family of eighth-order derivative-free methods. The proper selection  
of weight functions showed a reasonable reduction in error as compared to other referenced derivative-free methods.
It is obvious that constructed family has broad choice for the weight function in the third and fourth step of the method.
The true essence of the family is hidden in the construction of interpolation polynomial for the approximation of $f'(z)$
and weight functions make it more flexible to get higher performance and efficiency.

\bibliographystyle{model1-num-names}

\begin{thebibliography}{00}


\bibitem{1}
 J. F. steffensen, ''Remarks on iteration'', Skand Aktuar Tidsr, vol. 16, pp. 64-72, 1933.
\bibitem{2}
Y. khan, M. Faridi, K. Sayevand, A new general eighth-order family of iterative methods for solving nonlinear equations, Applied Mathematics 
Letters, 25 (2012) 2262-2266
\bibitem{3}
Jisheng Kou, Xiuhua Wang, Yitian Li, Some eighth-order root-finding three-step methods, Commun. Nonlinear Sci. Numer. Simul. 15 (2010) 536-544.
\bibitem{4}
S. Weerakoon and T. G. I. Fernando. A variant of newtonâs method with accelerated third-order convergence. Appl. Math.
Lett., 13:87:93, 2000.
\bibitem{5}
H. Kung and J. F. Traub. Optimal order of one-point and multipoint iteration. J. Assoc. Comput. Math., 21:643:651, 1974.
\bibitem{6}
J. F. Traub. Iterative Methods for solution of equations. Chelsea publishing company, New York, 1977.
\bibitem{7}
R. Thukral. A Family of Three-Point Derivative-Free Methods of Eighth-order for Solving Nonlinear Equations. J. of Mod. Meth. in Numer. Math., 
Vol. 3, No. 2, 2012, 11â21



 \end{thebibliography}

\end{document}